\documentclass{amsart}
\usepackage{amsmath,amsthm,amssymb,amsfonts,amscd, array}
\usepackage{epsfig}
\usepackage[all]{xy}
\newtheorem{lem}{Lemma}[section]
\newtheorem{prop}[lem]{Proposition}
\newtheorem{cor}[lem]{Corollary}
\newtheorem{thm}[lem]{Theorem}

\newcommand{\Ext}{\operatorname{Ext}\nolimits}

\newcommand{\mo}{\operatorname{mod}\nolimits}

\newcommand{\End}{\operatorname{End}\nolimits}

\newcommand{\op}{\operatorname{op}\nolimits}

\begin{document}
\title{Counting cluster-tilted algebras of type $A_n$}
\author{Hermund Andr\' e Torkildsen}
\address{Hermund Andr\' e Torkildsen, Department of Mathematical Sciences, Norwegian University of Science and
  Technology (NTNU), 7491 Trondheim, Norway}
\email{hermund.torkildsen@math.ntnu.no}

\begin{abstract}
The purpose of this paper is to give an explicit formula for the
number of non-isomorphic cluster-tilted algebras of type $A_n$, by
counting the mutation class of any quiver with underlying graph
$A_n$. It will also follow that if $T$ and $T'$ are cluster-tilting
objects in a cluster category $\mathcal{C}$, then $\End_{\mathcal{C}}(T)$ is
isomorphic to $\End_{\mathcal{C}}(T')$ if and only if $T=\tau^i
T'$. 
\end{abstract}

\maketitle

\section{Cluster-tilted algebras}
The cluster category was introduced independently in \cite{ccs} for
type $A_n$ and in \cite{bmrrt} for the general case. Let
$\mathcal{D}^b (\mo H)$ be the bounded derived category of the
finitely generated modules over a finite dimensional hereditary
algebra $H$ over a field $K$. In \cite{bmrrt} the cluster category was
defined as the orbit category $\mathcal{C}=\mathcal{D}^b (\mo H) /
\tau^{-1}[1]$, where $\tau$ is the Auslander-Reiten translation and
[1] the suspension functor. The cluster-tilted algebras are the
algebras of the form $\Gamma=\End_{\mathcal{C}}(T)^{\op}$, where $T$
is a cluster-tilting object in $\mathcal{C}$. See \cite{bmr1}. 
 
Let $Q$ be a quiver with no multiple arrows, no loops and no
oriented cycles of length two. Mutation of $Q$ at vertex $k$ is a
quiver $Q'$ obtained from $Q$ in the following way.

\begin{enumerate}
\item Add a vertex $k^{*}$.
\item If there is a path $i\rightarrow k \rightarrow j$, then if there
  is an arrow from $j$ to $i$, remove this arrow. If there is no arrow
  from $j$ to $i$, add an arrow from $i$ to $j$.
\item For any vertex $i$ replace all arrows from $i$ to $k$ with
  arrows from $k^{*}$ to $i$, and replace all arrows from $k$ to $i$
  with arrows from $i$ to $k^{*}$.
\item Remove the vertex $k$.
\end{enumerate}

We say that a quiver $Q$ is
mutation equivalent to $Q'$, if $Q'$ can be obtained from $Q$ by a
finite number of mutations. The mutation class of $Q$ is all quivers
mutation equivalent to $Q$. It is known from \cite{fz3} that the
mutation class of a Dynkin quiver $Q$ is finite.

If $\Gamma$ is a cluster-tilted algebra, then we say that $\Gamma$ is
of type $A_n$ if it arises from the cluster category of a path algebra
of Dynkin type $A_n$.  

Let $Q$ be a quiver of a cluster-tilted algebra $\Gamma$. From \cite{bmr2}, it
is known that if $Q'$ is obtained from $Q$ by a finite number of
mutations, then there is a cluster-tilted algebra $\Gamma '$ with quiver
$Q'$. Moreover, $\Gamma$ is of finite representation type if and only if
$\Gamma '$ is of finite representation type \cite{bmr1}. We also have that
$\Gamma$ is of type $A_n$ if and only if $\Gamma '$ is of type
$A_n$. From \cite{bmr3} we know that a cluster-tilted algebra is up
to isomorphism uniquely determined by its quiver. See also \cite{ccs2}.

It follows from this that to count the number of cluster-tilted
algebras of type $A_n$, it is enough to count the mutation class of
any quiver with underlying graph $A_n$. 

\section{Category of diagonals of a regular $n+3$ polygon} 
We recall some results from \cite{ccs}. 

Let $n$ be a positive
integer and let $\mathcal{P}_{n+3}$ be a regular polygon with $n+3$
vertices. A diagonal is a straight line between two non-adjacent
vertices on the border. A triangulation is a maximal set of
diagonals which do not cross. If $\Delta$ is any triangulation of
$\mathcal{P}_{n+3}$, we know that $\Delta$ consists of exactly $n$
diagonals.   

Let $\alpha$ be a diagonal between vertex $v_1$ and vertex $v_2$ on
the border of $\mathcal{P}_{n+3}$. In \cite{ccs} a
\textit{pivoting elementary move} $P(v_1)$ is an anticlockwise move of
$\alpha$ to another diagonal $\alpha'$ about $v_1$. The vertices of
$\alpha'$ are $v_1$ and $v_2'$, where $v_2$ and $v_2'$ are vertices of
a border edge and rotation is anticlockwise. A \textit{pivoting path} from $\alpha$ to
$\alpha'$ is a sequence of pivoting elementary moves starting at
$\alpha$ and ending at $\alpha'$. 

Fix a positive integer $n$. Categories of diagonals of regular
$(n+3)$-polygons were introduced in \cite{ccs}. Let $\mathcal{C}_n$ be
the category with indecomposable objects all diagonals of the
polygon, and we take as objects formal direct sums of these
diagonals. Morphisms from $\alpha$ to $\alpha'$ are generated by
elementary pivoting moves modulo the mesh relations, which are defined
as follows. Let $\alpha$ and $\beta$ be diagonals, with $a$ and $b$
the vertices of $\alpha$ and $c$ and $d$ the vertices of
$\beta$. Suppose $P(c)P(a)$ takes $\alpha$ to $\beta$. Then
$P(c)P(a)=P(d)P(b)$. Furthermore, if one of the intermediate edges in
a pivoting elementary move is a border edge, this move is zero. It is
shown in \cite{ccs} that this category is equivalent to the cluster
category defined in Section 1 in the $A_n$ case.  

We have the following from \cite{ccs}.

\begin{itemize}
\item The irreducible morphisms in $\mathcal{C}_n$ are the
  direct sums of pivoting elementary moves.
\item The Auslander-Reiten translation of a diagonal is given by
  clockwise rotation of the polygon. 
\item $\Ext_{\mathcal{C}_n}^{1}(\alpha,\alpha ') =
  \Ext_{\mathcal{C}}^{1}(\alpha,\alpha ') = 0$ if and only if
  $\alpha$ and $\alpha '$ do not cross.
\end{itemize}

It follows that a tilting object in $\mathcal{C}$ corresponds to a
triangulation of $\mathcal{P}_{n+3}$. 

For any triangulation $\Delta$ of $\mathcal{P}_{n+3}$, it is possible to define a quiver
$Q_{\Delta}$ with $n$ vertices in the following way. The
vertices of $Q_{\Delta}$ are the midpoints of the diagonals of
$\Delta$. There is an arrow between $i$ and $j$ in $Q_{\Delta}$ if the
corresponding diagonals bound a common triangle. The orientation is $i
\rightarrow j$ if the diagonal corresponding to $j$ is obtained from
the diagonal corresponding to $i$ by rotating anticlockwise about
their common vertex. It is known from \cite{ccs} that all quivers
obtained in this way are quivers of cluster-tilted algebras of type
$A_n$. 

We defined the  mutation of a quiver of a cluster-tilted algebra above. We
also define mutation of a triangulation at a given diagonal, by
replacing this diagonal with another one. This can be done in one and
only one way. Let $Q_{\Delta}$ be a quiver corresponding to a triangulation
$\Delta$. Then mutation of $Q_{\Delta}$ at the vertex $i$ corresponds
to mutation of $\Delta$ at the diagonal corresponding to $i$.

It follows that any triangulation gives rise to a quiver of a
cluster-tilted algebra, and that a quiver of a cluster-tilted algebra
can be associated to at least one triangulation.

Let $\mathcal{M}_n$ be the mutation class of $A_n$, i.e. all quivers
obtained by repeated mutation from $A_n$, up to isomorphisms of
quivers. Let $\mathcal{T}_n$ be the set of all triangulations of
$\mathcal{P}_{n+3}$. We can define a function $\gamma : \mathcal{T}_n
\rightarrow \mathcal{M}_n$, where we set $\gamma(\Delta)=Q_{\Delta}$
for any triangulation $\Delta$ in $\mathcal{T}_n$. Note that $\gamma$
is surjective.  

\section{Counting cluster-tilted algebras of type $A_n$} 
If $a$ and $b$ are vertices on the border of a regular polygon, we say
that the \textit{distance} between $a$ and $b$ is the smallest number of border
edges between them. Let us say that a diagonal from $a$ to $b$ is
\textit{close to the border} if the distance between $a$ and $b$ is
exactly $2$. For a quiver $Q_{\Delta}$ corresponding to a
triangulation $\Delta$, let us always write $v_\alpha$ for the vertex
of $Q_{\Delta}$ corresponding to the diagonal $\alpha$. 

If $Q$ is a quiver of a cluster-tilted algebra of type $A_n$, we we
have the following facts \cite{bv,ccs,s}.

\begin{itemize}
\item All cycles are oriented.
\item All cycles are of length 3.
\item There does not exist two cycles that share one arrow.
\end{itemize}

\begin{lem}\label{sinksourcecycle}
If a diagonal $\alpha$ of a triangulation $\Delta$ is close to the border, then
the corresponding vertex $v_{\alpha}$ in $\gamma(\Delta)=Q_{\Delta}$ is either a source, a
sink or lies on a cycle (oriented of length $3$). 
\end{lem}
\begin{proof}
All cycles are oriented and of length $3$ in the $A_n$
case. Suppose that $\alpha$ is a diagonal in $\Delta$ which is close
to the border. There are only three cases to consider, shown in Figure
\ref{figsinksourcecycle}.

  \begin{figure}[htp]
  \begin{center}
    \includegraphics[width=3.3cm]{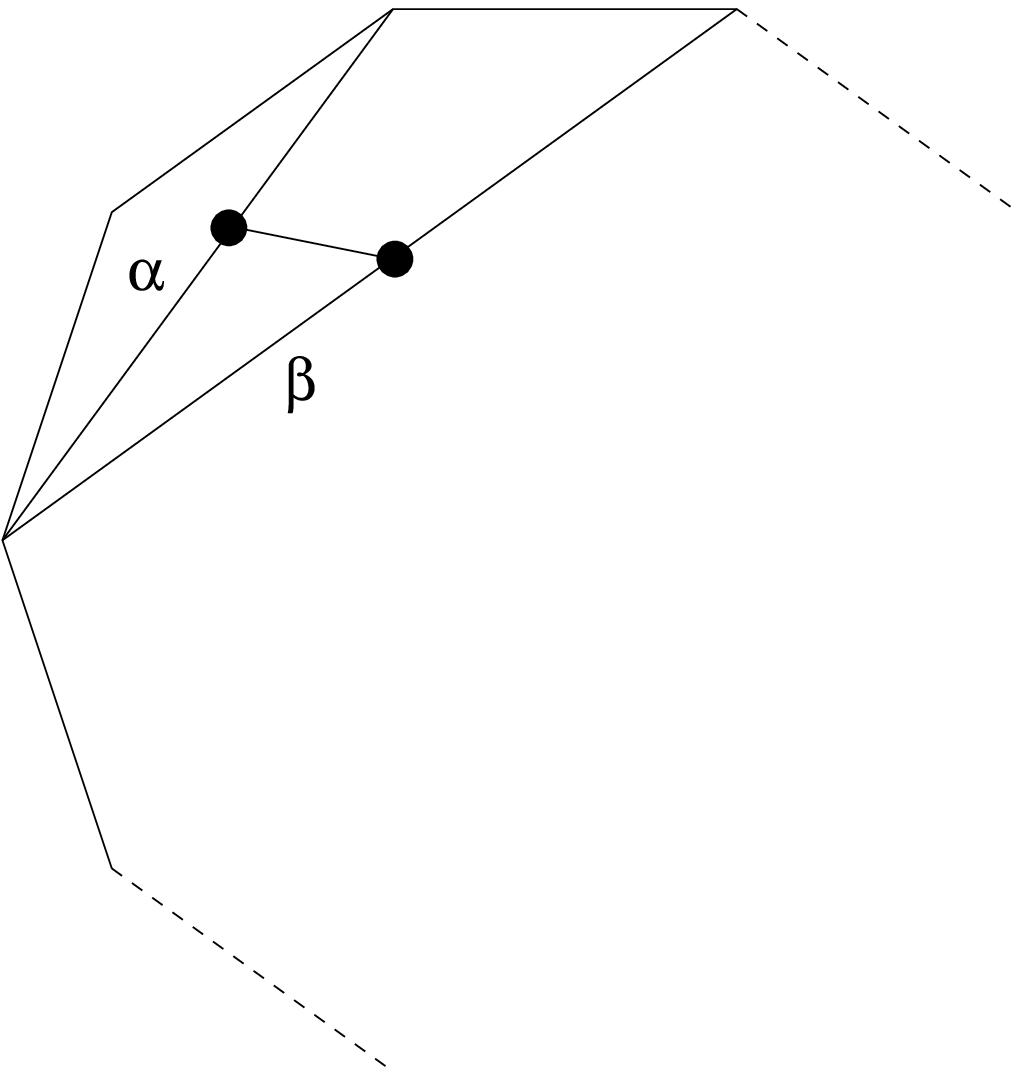}
    \includegraphics[width=3.3cm]{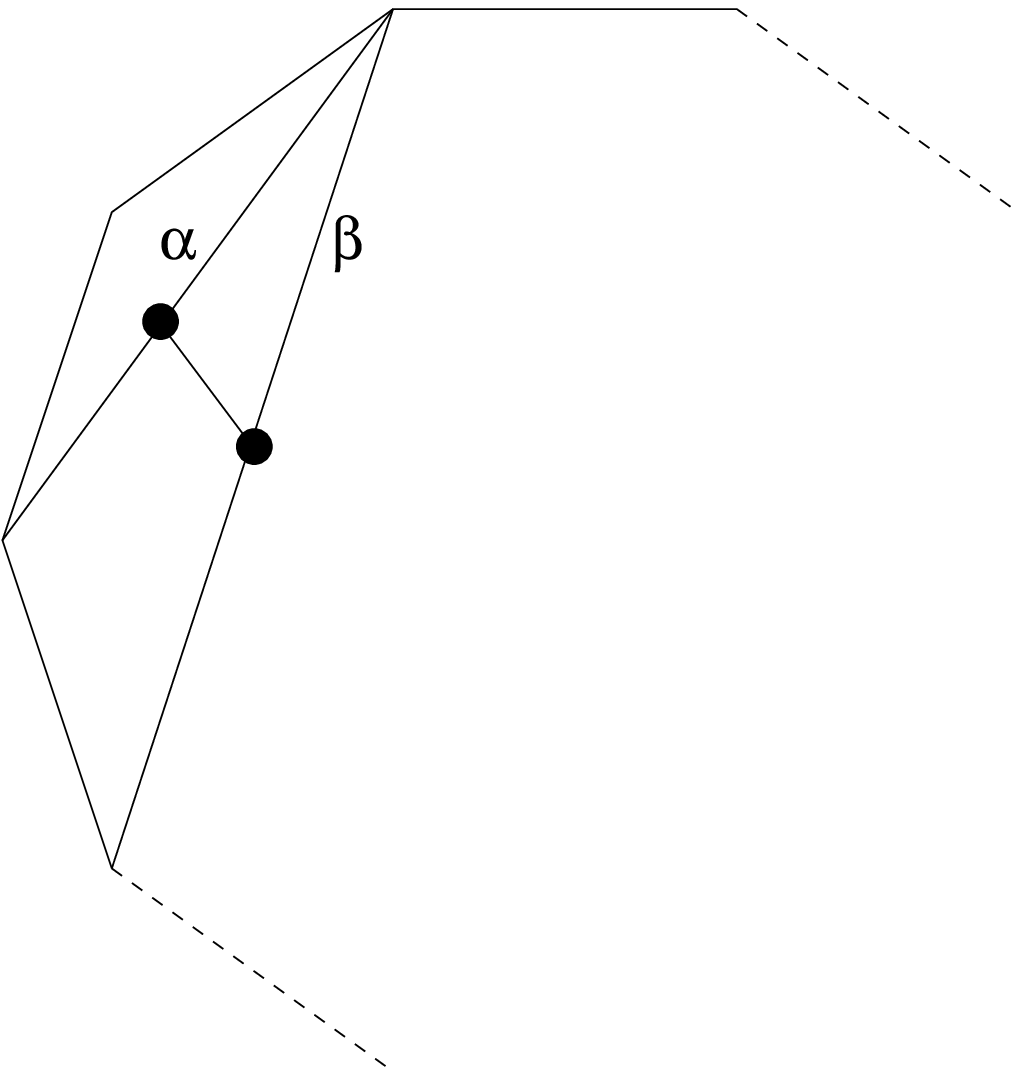}
    \includegraphics[width=3.3cm]{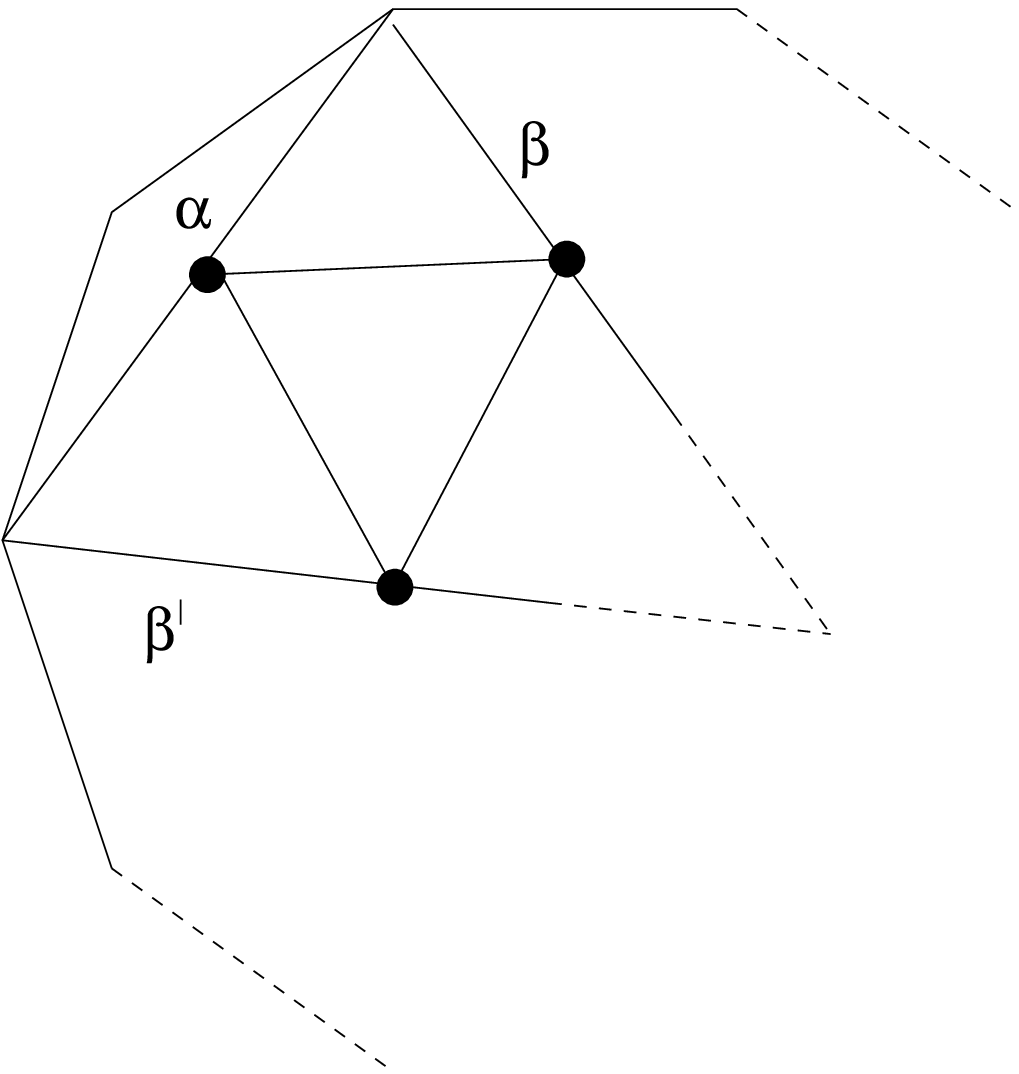}
  \end{center}\caption{\label{figsinksourcecycle}See the proof of Lemma \ref{sinksourcecycle}. Sink, source and cycle.}
  \end{figure}

  In the first case $\alpha$ corresponds to a sink. There is no other
  vertex adjacent to $v_{\alpha}$ but $v_{\beta}$, or else the
  corresponding diagonal of this vertex would cross $\beta$. We have
  the same for the second case where $\alpha$ is a source. In the
  third case $v_{\alpha}$ lies on a cycle.  
\end{proof}

Note that if $v_{\alpha}$ is a sink (or source) then $v_{\alpha}$ has
only one adjacent vertex if and only if $\alpha$ is close to the
border.
  
\begin{lem}\label{lemconnected}
Let $\Delta$ be a triangulation and let $\gamma(\Delta)=Q_{\Delta}$ be the
corresponding quiver. A quiver $Q'$ obtained from $Q_{\Delta}$ by factoring out
a vertex $v_{\alpha}$ is connected if and only if the corresponding diagonal 
$\alpha$ is close to the border.
\end{lem}
\begin{proof}
Suppose $\alpha$ is close to the border. By Lemma
\ref{sinksourcecycle}, $\alpha$ corresponds to a sink, a source or a
vertex on a cycle. If $v_\alpha$ is a sink or a source then $v_\alpha$
has only one adjacent vertex, so factoring out $v_\alpha$ does not
disconnect the quiver. Suppose $v_\alpha$ lies on a cycle. Then we are
in the case shown in the third picture in Figure
\ref{figsinksourcecycle}. We see that there can be no other vertex
adjacent to $v_\alpha$ except $v_\beta$ and $v_{\beta'}$, since else
the corresponding diagonal would cross $\beta$ or $\beta'$. Hence
factoring out $v_\alpha$ does not disconnect the quiver.
 
Next, suppose that factoring out $v_{\alpha}$ does not disconnect the
quiver. If $v_{\alpha}$ is a source or a sink with only one adjacent
vertex, then $v_{\alpha}$ is close to the border. If not, first suppose
$v_{\alpha}$ does not lie on a cycle. Then it is clear that factoring
out $v_{\alpha}$ disconnects the quiver, so we may assume that
$v_{\alpha}$ lies on a cycle. Then $\alpha$ is an edge of a triangle
consisting of only diagonals (i.e. no border edges), say $\beta$ and
$\beta'$. Suppose there is a vertex $v_{\delta}$ adjacent to $v_{\alpha}$,
with $v_{\delta} \neq v_{\beta}$ and $v_{\delta} \neq v_{\beta'}$. Then $v_{\delta}$ can not be
adjacent to $v_{\beta}$ or $v_{\beta'}$, since then we would have two
cycles sharing one arrow. We also see that $v_{\delta}$ can not be adjacent to
any vertex $v_{\gamma}$ from which there exists a path to
$v_{\beta}$ or $v_{\beta'}$ not containing $v_{\alpha}$, or else there
would be a cycle of length greater than $3$. Therefore factoring out
$v_{\alpha}$ would disconnect the quiver, and this is a contradiction,
thus there can be no other vertices adjacent to $v_{\alpha}$. It
follows that $\alpha$ can not be adjacent to any other diagonal but
$\beta$ and $\beta'$, hence $\alpha$ is close to the border.  
\end{proof}

\begin{figure}[h]
  \begin{center}
    \includegraphics[width=5cm]{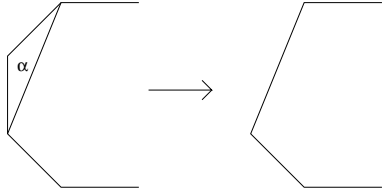}
  \end{center}\caption{\label{figfactoring}Factoring
    out a diagonal close to the border}
  \end{figure}

Let $\Delta$ be a triangulation of $\mathcal{P}_{n+3}$ and let
$\alpha$ be a diagonal close to the border. The triangulation
$\Delta'$ of $\mathcal{P}_{n+3-1}$ obtained from $\Delta$ by factoring
out $\alpha$ is defined as the triangulation of $\mathcal{P}_{n+3-1}$
by letting $\alpha$ be a border edge and leaving all the other
diagonals unchanged. We write $\Delta / \alpha$ for the new
triangulation obtained. See Figure \ref{figfactoring}.  

\begin{lem}\label{lemfactoring}
Let $\Delta$ be a triangulation and $\gamma(\Delta)=Q_{\Delta}$.
Factoring out a vertex in $Q_{\Delta}$ such that the resulting quiver
is connected, corresponds to factoring out a diagonal of $\Delta$
close to the border.    
\end{lem}
\begin{proof}
Factoring out a vertex $v_{\alpha}$ in $Q$ such that the resulting quiver is
connected, implies that $\alpha$ is close to the border by Lemma
\ref{lemconnected}. Then consider all cases shown in Figure
\ref{figsinksourcecycle}.   
\end{proof}

Note that this means that $\gamma(\Delta / \alpha) = Q_{\Delta} /
v_{\alpha}$. We have the following easy fact.

\begin{prop}\label{factor out}
Let $Q$ be a quiver of a cluster-tilted algebra of type $A_n$, with $n
\geq 3$. Let $Q'$ be obtained from $Q$ by factoring out a vertex such
that $Q'$ is connected. Then $Q'$ is the quiver of some cluster-tilted
algebra of type $A_{n-1}$.    
\end{prop}
\begin{proof}
It is already known from \cite{bmr2} that $Q'$ is the quiver of a
cluster-tilted algebra. Suppose $\Delta$ is a triangulation of
$\mathcal{P}_{n+3}$ such that $\gamma(\Delta)=Q$. Such a $\Delta$ exists since
$\gamma$ is surjective. It is enough, by Lemma
\ref{lemconnected}, to consider vertices corresponding to a diagonal
close to the border. By Lemma \ref{lemfactoring}, factoring out a
vertex corresponding to a diagonal $\alpha$ close to the border,
corresponds to factoring out $\alpha$. Then the resulting
triangulation of $\mathcal{P}_{(n-1)+3}$ corresponds to a quiver of a
cluster-tilted algebra of type $A_{n-1}$, since it is a triangulation.
\end{proof}

\begin{figure}[htp]
  \begin{center}
    \includegraphics[height=3.3cm]{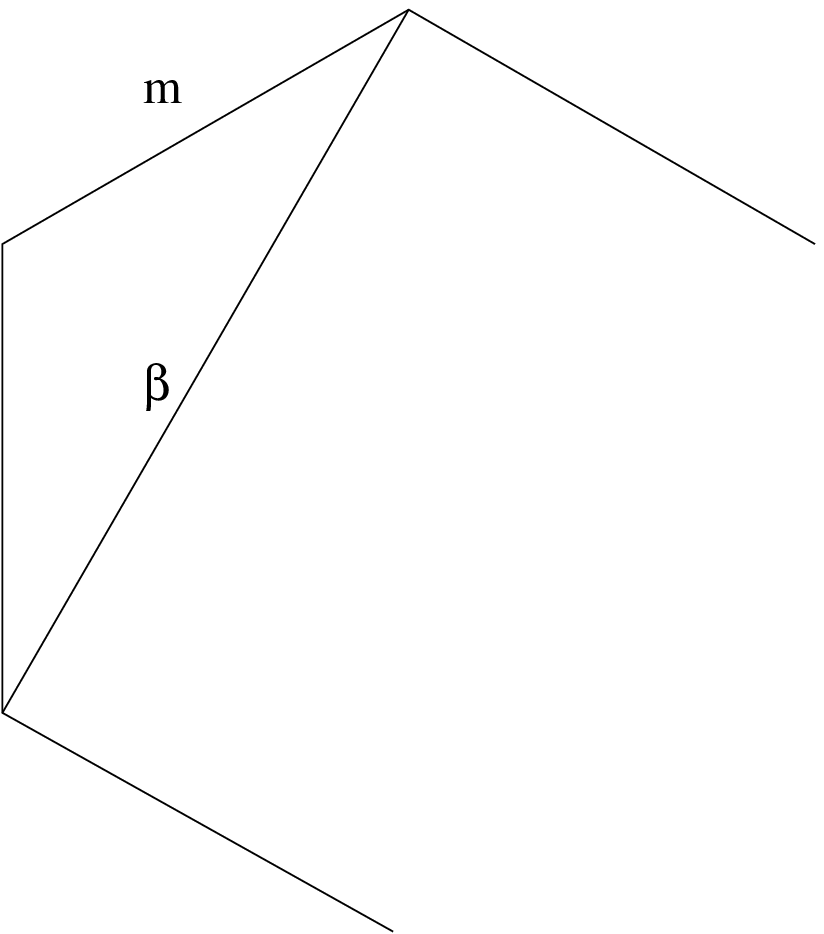}
    \includegraphics[height=3.3cm]{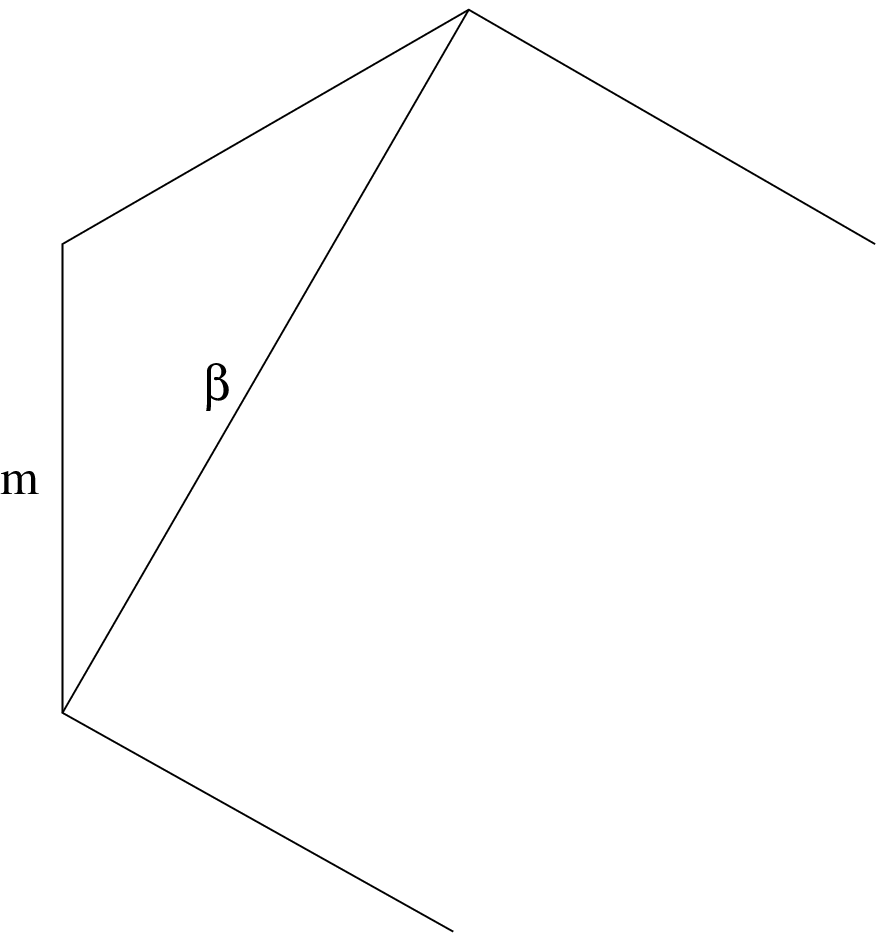}
    \includegraphics[height=3.3cm]{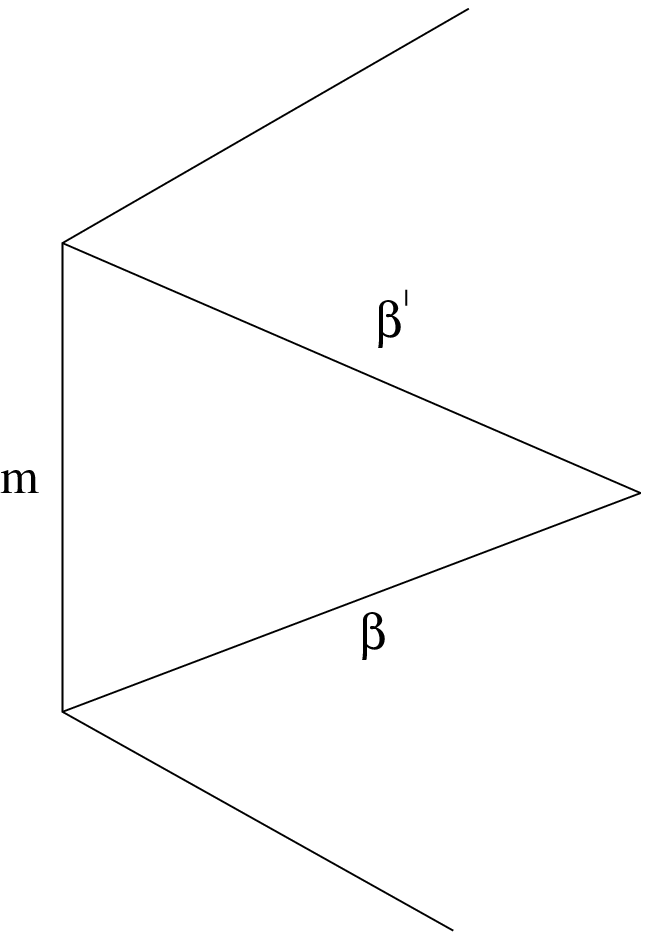}
  \end{center}\caption{\label{figBorderCases}}
  \end{figure}

Now we want to do the opposite of factoring out a vertex close to the
border. If $\Delta$ is a triangulation of $\mathcal{P}_{n+3}$, we want
to add a diagonal $\alpha$ such that $\alpha$ is a diagonal close to
the border and such that $\Delta \cup \alpha$ is a triangulation of
$\mathcal{P}_{(n+1)+3}$. Consider any border edge $m$ on
$\mathcal{P}_{n+3}$. Then we have one of the cases shown in Figure
\ref{figBorderCases}.

\begin{figure}[htp]
  \begin{center}
    \includegraphics[height=3.3cm]{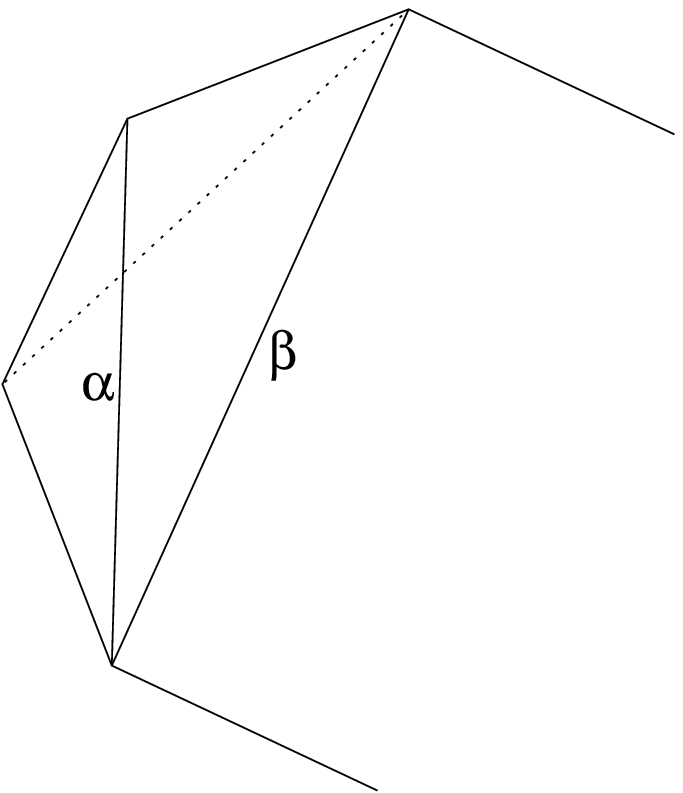}
    \includegraphics[height=3.3cm]{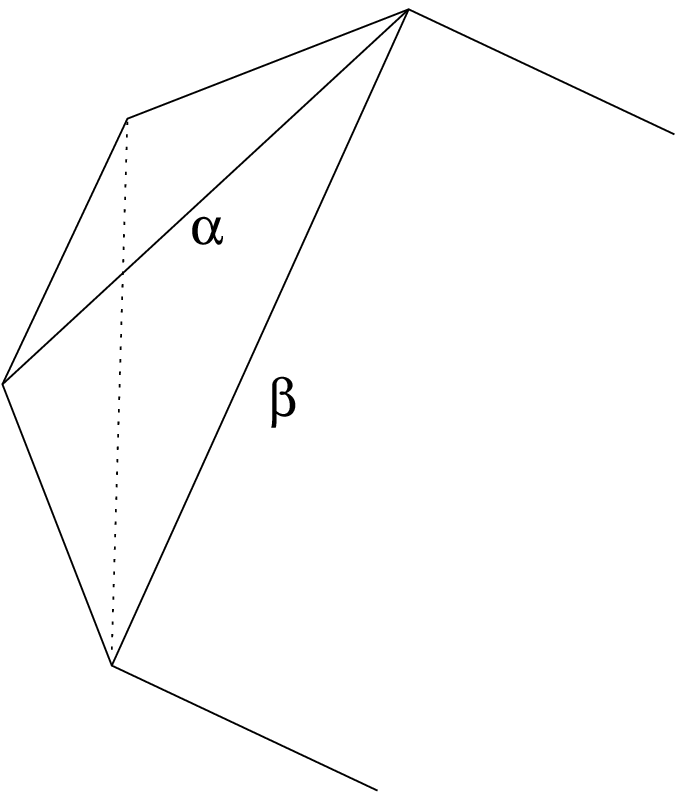}
    \includegraphics[height=3.3cm]{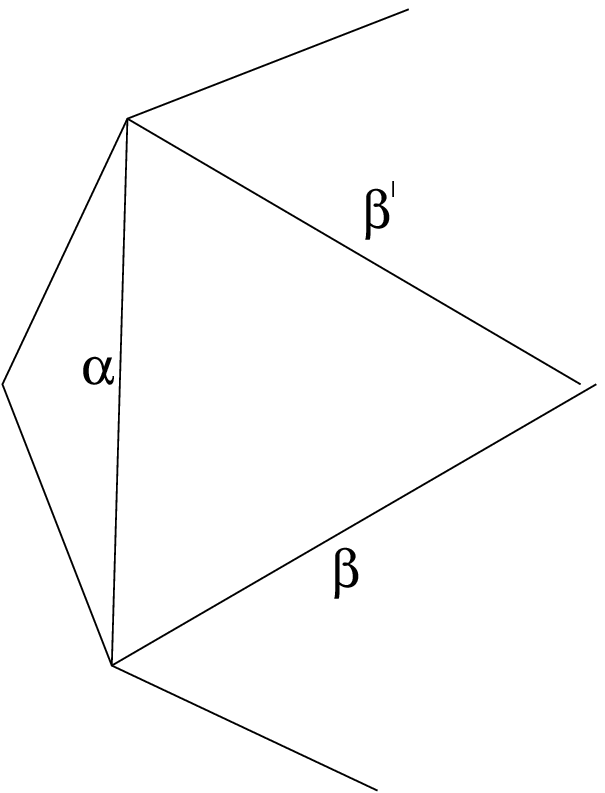}
  \end{center}\caption{\label{figExtension}}
  \end{figure}

We can extend the polygon at $m$ for each case in Figure
\ref{figBorderCases}, and add a diagonal $\alpha$ to the
extension. See Figure \ref{figExtension}
for the corresponding extensions at $m$.

It follows that for a given diagonal $\beta$, there are at
most three ways to extend the polygon with a diagonal $\alpha$ such
that $\alpha$ is adjacent to $\beta$, and it is easy to see that these
extensions gives non-isomorphic quivers.

For a triangulation $\Delta$ of $\mathcal{P}_{n+3}$, let us denote by
$\Delta(i)$ the triangulation obtained from $\Delta$ by rotating
$\Delta$ $i$ steps in the clockwise direction. We define an
equivalence relation on $\mathcal{T}_n$, where we let $\Delta
\sim \Delta(i)$ for all $i$. We define a new function
$\tilde{\gamma}: (\mathcal{T}_n / \sim) \rightarrow
\mathcal{M}_n$ induced from $\gamma$. This is well defined, for if
$\Delta = \Delta'(i)$ for an $i$, then obviously $Q_{\Delta} =
Q_{\Delta'}$ in $\mathcal{M}_n$. And hence since $\gamma$ is a
surjection, we also have that $\tilde{\gamma}$ is a surjection. We
actually have the following. 

\begin{thm}\label{quivers not iso}
The function $\tilde{\gamma}:(\mathcal{T}_n / \sim) \rightarrow
\mathcal{M}_n$ is bijective for all $n \ge 2$.
\end{thm}
\begin{proof}
We already know that $\tilde{\gamma}$ is surjective.

Suppose $\tilde\gamma(\Delta) = \tilde\gamma(\Delta')$ in $\mathcal{M}_n$. We
want to show that $\Delta = \Delta'$ in $(\mathcal{T}_n / \sim)$ using
induction.  

It is easy to check that $(\mathcal{T}_3 / \sim) \rightarrow
\mathcal{M}_3$ is injective. Suppose $(\mathcal{T}_{n-1} / \sim) \rightarrow
\mathcal{M}_{n-1}$ is injective. Let $\alpha$ be a diagonal close to the
border in $\Delta$, with image $v_{\alpha}$ in $Q$, where $Q$ is a
representative for $\tilde{\gamma}(\Delta)$. Then the
diagonal $\alpha '$ in $\Delta '$ corresponding to $v_{\alpha}$ in $Q$
is also close to the border. We have $\tilde{\gamma}(\Delta / \alpha) =
\tilde{\gamma}(\Delta' / \alpha ') = Q / v_{\alpha}$ by Lemma
\ref{lemfactoring}, and hence, by hypothesis, $\Delta / \alpha =
\Delta' / \alpha '$ in $(\mathcal{T}_n / \sim)$.

We can obtain $\Delta$ and $\Delta '$ from  $\Delta / \alpha =
\Delta' / \alpha '$
by extending the polygon at some border edge. Fix a diagonal $\beta$
in $\Delta$ such that
$v_{\alpha}$ and $v_{\beta}$ are adjacent. This can be done since $Q$
is connected. Let $\beta '$ be the diagonal in $\Delta'$ corresponding
to $v_{\beta}$. By the above there are at most three ways to extend
$\Delta / \alpha$ such that the new diagonal is adjacent to $\beta$. It is clear that these extensions will be mapped
by $\tilde{\gamma}$ to non-isomorphic quivers. Also there are at most
three ways to extend $\Delta' / \alpha'$ such that the new diagonal is
adjacent to $\beta'$, and all these extensions are mapped to
non-isomorphic quivers, thus $\Delta = \Delta'$ in $(\mathcal{T}_n / \sim)$.
\end{proof} 

Note that this also means that $\Delta=\Delta'(i)$ for an $i$
if and only if $Q_{\Delta} \simeq Q_{\Delta'}$ as quivers.

Now, let $T$ be a cluster-tilting object of the cluster category $\mathcal{C}$. This object
corresponds to a triangulation $\Delta$ of $\mathcal{P}_{n+3}$, and
all tilting objects obtained from rotation of $\Delta$ gives the same
cluster-tilted algebra. No other triangulation gives rise to the same
cluster-tilted algebra. 

The Catalan number $C(i)$ can be defined as the number of
triangulations of an $i$-polygon with $i-3$ diagonals. The number is
given by the following formula.

$$C(i)=\frac{(2i)!}{(i+1)!i!}$$ 

We now have the following.

\begin{cor}
The number $a(n)$ of non-isomorphic basic cluster-tilted algebras of type $A_n$ is
the number of triangulations of the disk with $n$ diagonals, i.e. 
$$a(n)=C(n+1)/(n+3)+C((n+1)/2)/2+(2/3)C(n/3),$$
where $C(i)$ is the $i$'th Catalan number and the second term is
omitted if $(n+1)/2$ is not an integer and the third term is omitted
if $n/3$ is not an integer. 
\end{cor}

These numbers appeared in a paper by W. G. Brown in 1964 \cite{b}. See
Table \ref{examples} for some values of $a(n)$.

\begin{table}
\begin{tabular}{l|l}
$n$&$a(n)$\\ \hline
2&1\\
3&4\\
4&6\\
5&19\\
6&49\\
\end{tabular}
\; \;
\begin{tabular}{l|l}
$n$&$a(n)$\\ \hline
7&150\\
8&442\\
9&1424\\
10&4522\\
11&14924\\
\end{tabular}
\caption{Some values of $a(n)$.}\label{examples}
\end{table}
 
We have that if $T$ is a cluster-tilting object in $\mathcal{C}$, then the
cluster-tilted algebras $\End_{\mathcal{C}}(T)$ and
$\End_{\mathcal{C}}(\tau T)$ are isomorphic. In the $A_n$ case we also
have the following.

\begin{thm}\label{ct not iso}
Let $T$ and $T'$ be tilting objects in $\mathcal{C}$, then the
cluster-tilted algebras $\End_{\mathcal{C}}(T)$ and
$\End_{\mathcal{C}}(T')$ are isomorphic if and only if $T' =
\tau^{i}T$ for an $i \in \mathbb{Z}$.
\end{thm} 
\begin{proof}
Let $\Delta$ be the triangulation of $\mathcal{P}_{n+3}$ corresponding
to $T$ and let $\Delta'$ be the triangulation corresponding to
$T'$. If $T' \not\simeq \tau^{i}T$ for any $i$, then $\Delta'$ is not
obtained from $\Delta$ by a rotation, and hence
$\End_{\mathcal{C}}(T)$ is not isomorphic to $\End_{\mathcal{C}}(T')$
by Theorem \ref{quivers not iso}. 
\end{proof} 

\begin{prop}
Let $\Gamma$ be a cluster-tilted algebra of type $A_n$. The number of
non-isomorphic cluster-tilting objects $T$ such that $\Gamma \simeq
\End_{\mathcal{C}}(T)$ has to divide $n+3$.  
\end{prop}
\begin{proof}
Let $T$ be a tilting object in $\mathcal{C}$ corresponding to the
triangulation $\Delta$. Denote by $\Delta(i)$ the rotation of $\Delta$
$i$ steps in the clockwise direction. Let $0 < s \leq n$ be the smallest number of
rotations needed to obtain the same triangulation $\Delta$, i.e. the
smallest $s$ such that $\Delta = \Delta (s)$. It is clear from the above that $T \not\simeq
T'$, where $T'$ corresponds to $\Delta (t)$  with $0 < t < s$, hence
$s$ is the number of non-isomorphic tilting objects giving the same cluster-tilted
algebra. Now we only need to show that $s$ divides $n+3$, but this is
clear.
\end{proof}

The proof of the following is easy and is left to the reader. First
recall from \cite[Proposition 3.8]{fz2} that there are exactly $C(n)$ non-isomorphic
tilting objects in the cluster category for type $A_n$, where $C(n)$
denotes the $n$'th Catalan number.

\begin{prop} Consider the $A_n$ case.
\begin{itemize}
\item There are always at least $2$ non-isomorphic cluster-tilting objects giving the same
  cluster-tilted algebra.
\item There are at most $n+3$ non-isomorphic cluster-tilting objects giving the same
  cluster-tilted algebra.
\item Let $\Gamma$ be a cluster-tilted algebra of type $A_n$. If $n+3$ is prime,
  there are exactly $n+3$ non-isomorphic cluster-tilting objects giving $\Gamma$. In
  this case there are $C(n)/n+3$ non-isomorphic cluster-tilted
  algebras, where $C(n)$ denotes the $n$'th Catalan number.
\end{itemize}
\end{prop}
\begin{flushleft}
\textbf{Acknowledgements:} I would like to thank the referee and my supervisor
Aslak Bakke Buan for valuable comments. 
\end{flushleft}

\end{document}